\documentclass{elsarticle}

\usepackage{amsmath, amssymb}

\newdefinition{defn}{Definition}

\newcommand{\R}{\mathbb{R}}
\newcommand{\set}[1]{\begin{Bmatrix}{#1}\end{Bmatrix}}
\renewcommand{\r}{\mathcal{R}}

\DeclareMathOperator{\PE}{Percent Error}

\begin{document}

\begin{frontmatter}
\title{An algorithm for autonomously plotting solution sets in the presence of turning points}
\author[McG1]{Steven Pollack \corref{cor1}}
\ead{steven.pollack@mail.mcgill.ca}

\author[UoT]{Daniel Badali\corref{cor2}}
\ead{daniel.badali@utoronto.ca}

\author[McG2]{Jonathan Pollack}
\ead{jon.pollack@mail.mcgill.ca}

\cortext[cor1]{Corresponding author}

\cortext[cor2]{Principal corresponding author}

\address[McG1]{Department of Mathematics, McGill University, Montreal, Quebec, Canada, H3A 2T5}

\address[UoT]{Department of Chemical and Physical Science, University of Toronto at Mississauga, Mississauga, Ontario, Canada, L5L 1C6}

\address[McG2]{Department of Mechanical Engineering, McGill University, Montreal, Quebec, Canada, H3A 2T5}

\date{}

\begin{abstract}
Plotting solution sets for particular equations may be complicated by the existence of turning points. Here we describe an algorithm which not only overcomes such problematic points, but does so in the most general of settings. Applications of the algorithm are highlighted through two examples: the first provides verification, while the second demonstrates a non-trivial application. The latter is followed by a thorough run-time analysis. While both examples deal with bivariate equations, it is discussed how the algorithm may be generalized for space curves in $\R^{3}$.
\end{abstract}


\end{frontmatter}

\section{Introduction}

In this paper we consider curves determined by equations of the form
\begin{equation}
f(x,y)=0
\label{eq:main}
\end{equation}
where $f : I \to \R$, $I \subset \R \times \R$ is a product of open intervals. Equations such as \eqref{eq:main} are often referred to as implicit equations since in general there does not exist an explicit, unique function $g$ such that $y = g(x)$. A canonical example is that of the unit circle, with $f(x,y) = x^2+y^2-1=0$, which cannot be rearranged to isolate $y$ as a function of $x$. The following discussion also applies to equations of the form $f(x;\alpha)=0$ where $\alpha$ is a bifurcation parameter.

The purpose of this paper is to address the problem of plotting solution curves to \eqref{eq:main} with turning points (a precise definition will be presented shortly). While a method for dealing with this problem has already been developed by Keller \cite{keller1,keller2} (via a pseudo-arc-length parametrization), it relies on the turning points to be anything but cusps, and thus cannot be used on curves without some understanding of their profile. Conversely, the algorithm that we will present not only requires no prior knowledge of the solution curve, but also approaches the problem in what we consider to be a more intuitive manner. 

It should be said, however, that although we acknowledge that there are other methods to deal with this problem, we will make no attempt to compare them. Essentially, this paper is to be self-contained, and therefore our only concern is the explanation of our proposed algorithm, and its own benefits, not its relative benefits.

\section{Background}
Before a discussion of the algorithm may begin, we must first define the term \textit{turning point}.

\begin{defn} Let $f : \R^2 \to \R$, and $S = \left\{ (x,y) \in \R^2 : f(x,y) = 0 \right\}$. A point $(x^*,y^*) \in S$ is a \textit{turning point} of $f(x,y)=0$ if $(x^*,y^*)\in S$, and there exists $\delta >0$ such that at least one of the following statements is true:

\begin{enumerate}[Type 1:] 
\item $\left\{(x,y) \in \R^2 : 0 < x - x^* < \delta, \; 0 < |y - y^*| < \delta\right\} \cap S = \varnothing$. 
\item  $\left\{(x,y) \in \R^2 : 0 < |x - x^*| < \delta, \; 0 < y - y^* < \delta\right\} \cap S = \varnothing$. 
\item  $\left\{(x,y) \in \R^2 : 0 < x^* - x < \delta, \; 0 < |y - y^*| < \delta\right\} \cap S = \varnothing$. 
\item  $\left\{(x,y) \in \R^2 : 0 < |x - x^*| < \delta, \; 0 < y^* - y < \delta\right\} \cap S = \varnothing$. 
\end{enumerate}
\end{defn}

\begin{figure}[ht!]
\centering 
\includegraphics[width=0.5\textwidth]{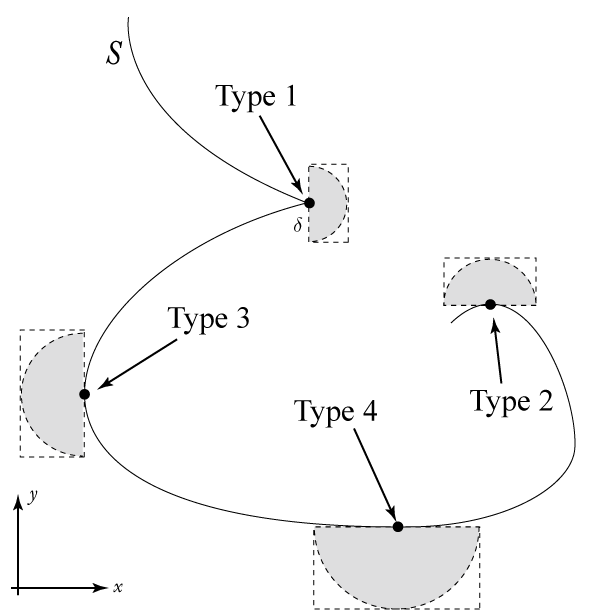}
\caption{Illustration of the different types of turning points}
\label{fig:tps}
\end{figure}

Definition 1 may be envisioned in the following way: if one were to restrict themselves to a small enough neighborhood about $(x^*, y^*)$, one of the open ``half-balls" centered at $(x^*, y^*)$ contained in this neighborhood would have no intersection with $S$ (see Figure \ref{fig:tps}). 

The remainder of this paper will deal with the discrete set:
\begin{equation}
S_k=\set{ \set{(x_j,y_j)}_{j=1}^{N} \subset \R^2 : f_k(x_j,y_j)=0, \; 1 \leq j \leq N }
\end{equation}
where $S_k$ is the approximation of $S$, and $f_k$ is the discrete approximation $f$. For a detailed discussion of approximating turning points in discrete spaces see the review by Cliffe, Spence, and Tavener\cite{cliffe}.

\section{Algorithm}

\begin{figure}
\centering
\includegraphics[scale=0.25]{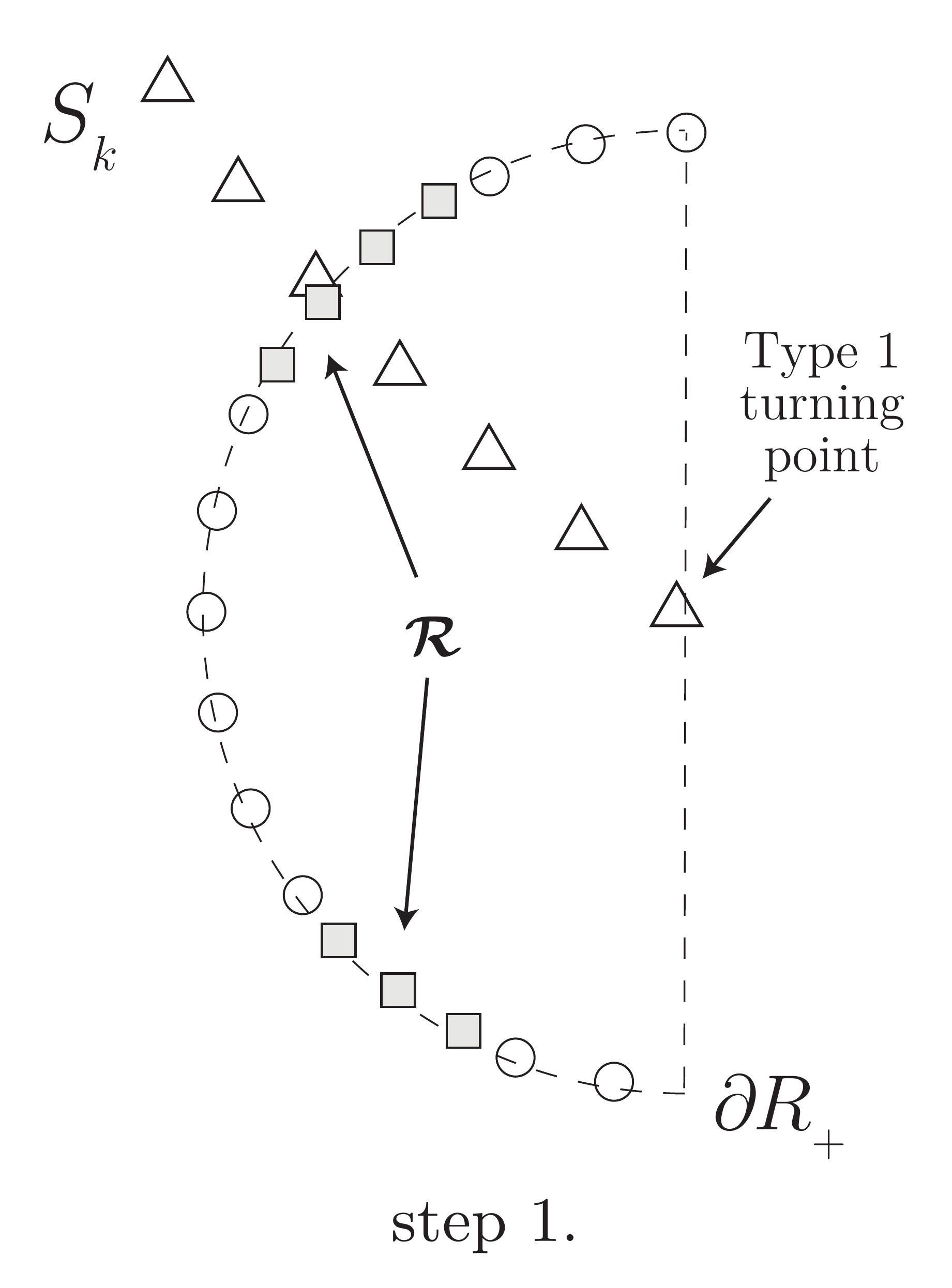}
\includegraphics[scale=0.25]{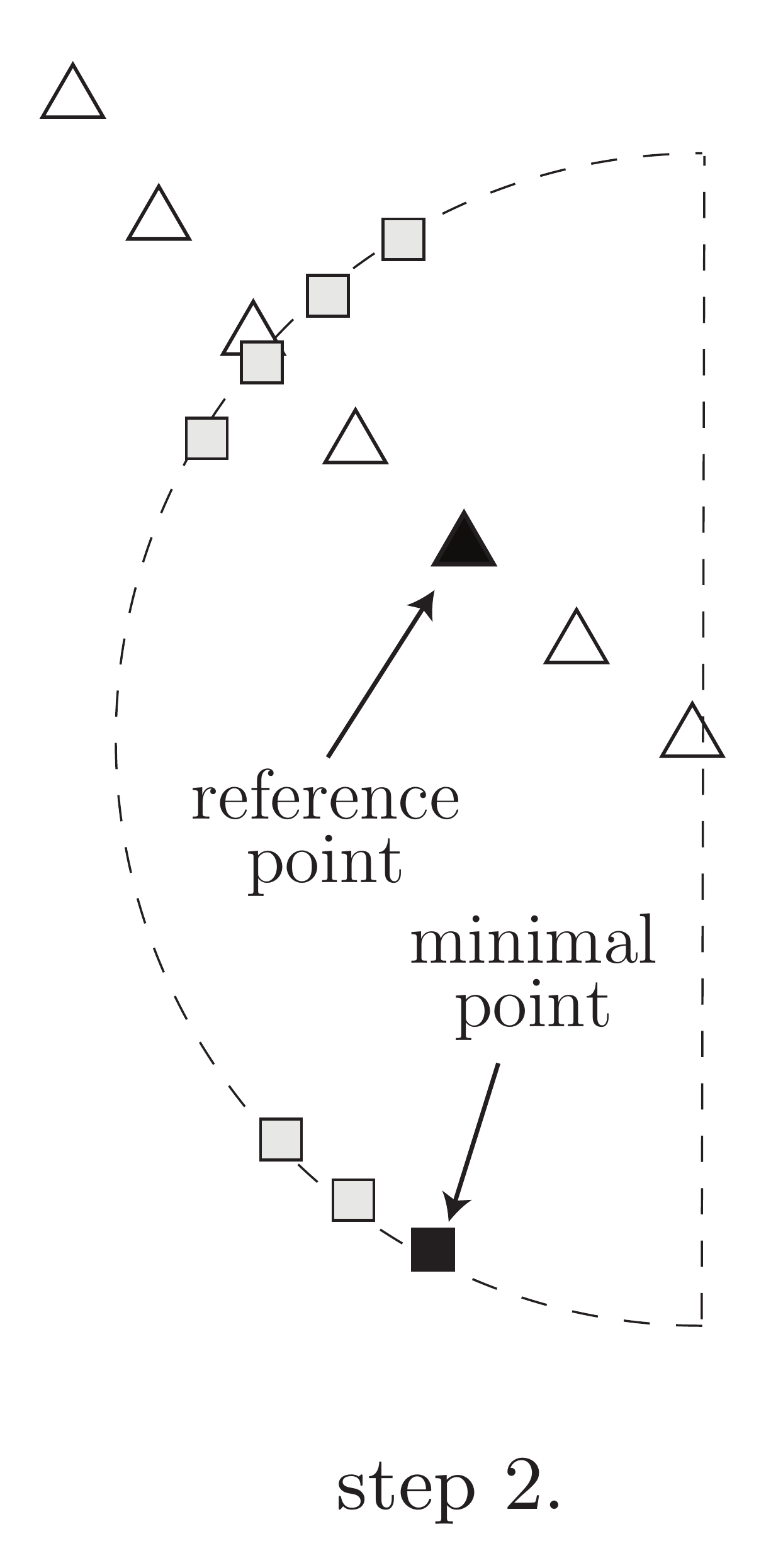}
\includegraphics[scale=0.25]{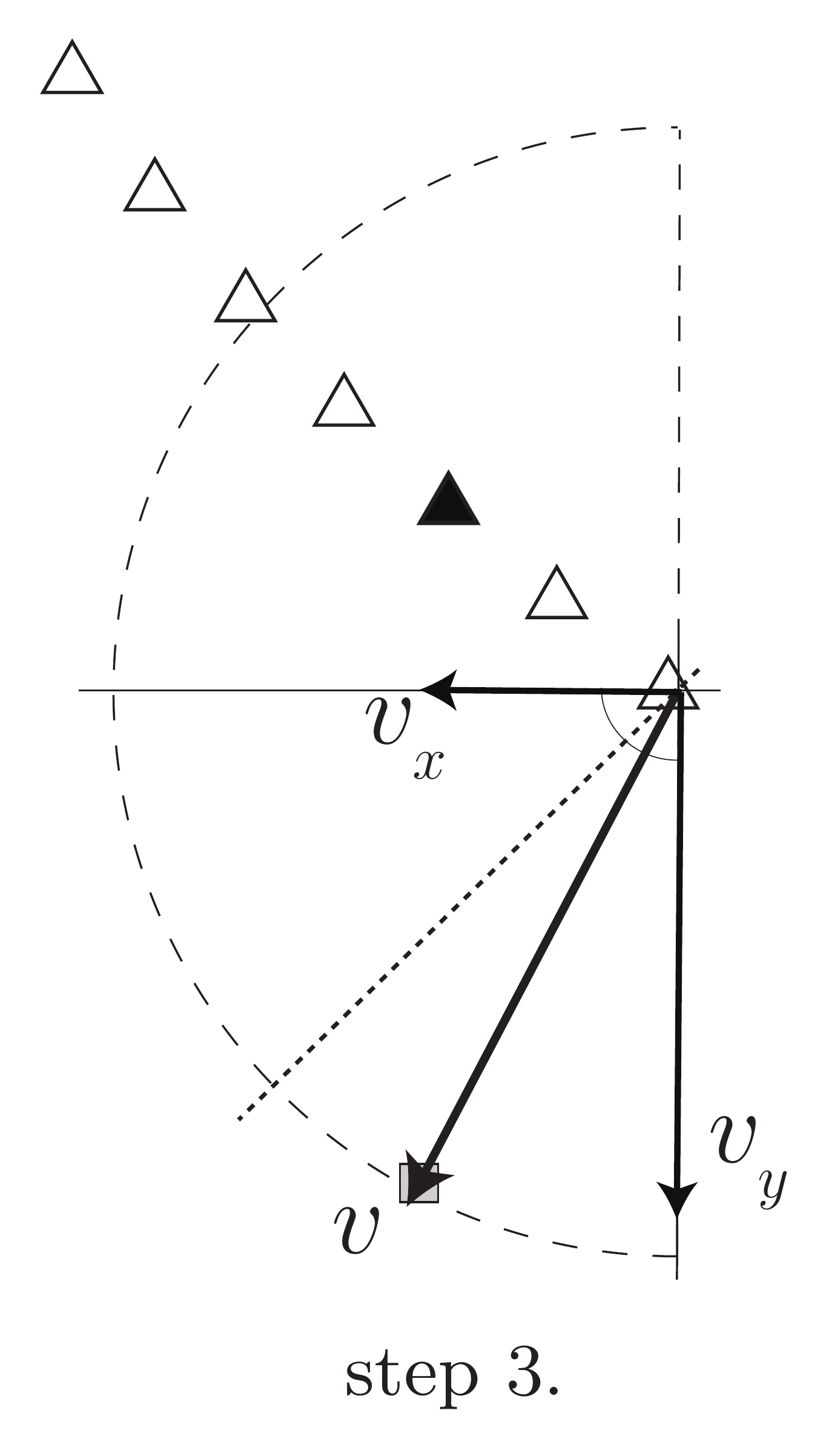}
\caption{Demonstration of algorithm on a Type 1 turning point.}
\label{fig:algo}
\end{figure}

Supposing that $(x_{j},y_{j}) \in S_k$ has been identified as a turning point of any type, let $R$ be the closed disk of radius $r$ centered at $(x_j, y_j)$. The algorithm is then performed in the following steps:
\begin{enumerate}[Step 1.]
\item Uniformly scan the boundary of $R$, $\partial R$, for solutions to $f(x,y) = 0$. Denote the set of solutions as $\r$. 
\item 
\begin{enumerate} \item If $\r = \varnothing$, we may stop here: no other points in $\R^2$ satisfy $f(x,y) = 0$.
\item  If $\r \neq \varnothing$, select $(x_i, y_i) \in S_k$, $i < j$, to be a ``reference point'' and select the minimal point $(s',t') \in \r$ such that $\phi(s',t') \leq \phi(s,t)$ for all $(s,t) \in \r$, where
\begin{equation}\label{cost_function}
\phi(s,t) = \frac{1}{\sqrt{(x_i - s)^2 + (y_i - t)^2}}
\end{equation}
\end{enumerate}
\item  
\begin{enumerate}
\item Create the vector $v = (s'-x_j, t'-y_j)$,
\item Determine the largest axial component of $v$, $d = \max\set{ |s'-x_j|, |t'-y_j|}$.
\item Use $d$ to determine the new direction of iteration.
\end{enumerate}
\end{enumerate}

To clear up any possible ambiguity, we will elaborate on how one might implement the above strategy on a Type 1 turning point.

\paragraph{Step 1} Since $(x_j, y_j)$ is a Type 1 turning point, we know that there exists an open ball of radius $\delta$ centered at $(x_j, y_j)$ whose eastern hemisphere is disjoint from $S_k$. Therefore, if any solution to the equation $f(x,y) = 0$ exists, it must exist in the western hemisphere. Hence, we describe $\partial R_{+}$, our reduced search path, by 
\[
\partial R_{+} = \set{ (x,y) \in \R^2 : (x-x_j)^2 + (y-y_j)^2 = \delta^2, \; x < x_j}
\]
which is easily seen to be parametrized by 
\[
g(\theta) = (\delta \cos(\theta)+x_j, \delta \sin(\theta) + y_j), \qquad \theta \in [\pi/2, 3\pi/2]
\]
Thus, if $n$ points are used to uniformly sample $\partial R_{+}$, the mesh would look like 
\[
g_i = \left(\delta \cos\left(\frac{\pi}{2n}(n+2i)\right) + x_j , \delta\sin\left(\frac{\pi}{2n}(n+2i)\right)+y_j \right) 
\]
for $i = 0, 1, \ldots, n-1$.

The points $g_i$ are then checked to see if $f(g_i)=0$. If $g_i$ is solution to $f(x,y) = 0$, it should be recorded in the set $\r$. It should be noted that the method employed to check whether $g_i$ is a solution to $f(x,y) = 0$ is irrelevant. 

\paragraph*{Step 2}  We have found that it is best to pick an $(x_i, y_i)$ that is contained in $R$. In this way, when we minimize the cost function $\phi$, see \eqref{cost_function}, we are searching for the the point in $\r$ which is furthest from our the reference point. If we assume $(x_j, y_j)$ is not the final point of $S_k$, then the furthest point from the reference point is not expected to be any of the previous points in $S_k$. Picking a reference point that is ``too far'' from $(x_j, y_j)$ can potentially lead to choosing $(s', t') \in \r$ such that $(s',t')$ is a previous solution. When the algorithm is implemented inside an iterated plotting program (one would do this for autonomous plotting), this has the consequence of leading the plotting program to retrace its old steps. 
 
 From here it is only a matter of evaluating $\phi(s,t)$ for all $(s,t) \in \r$, and finding the smallest value. The coordinate which minimizes $\phi$ shall be called $(s', t')$. While there is no garentee that there exists a single coordinate, such that $\phi$ is minimized uniquely; in practice, we seldom arrive at two, or more points, that equally minimize $\phi$.

\paragraph*{Step 3} Make the direction vector, $v$, from $(x_j, y_j)$ to $(s',t')$. Then identify which component of $v$ is larger in magnitude. If the $y$-component is larger, begin iterating along the $y$-axis (the direction is determined by the sign of the $y$-component). 

Iteration should thence begin at $(s',t')$. This is to compensate for the fact that it may happen that the turning point was approached via iteration along the positive $x$-direction, yet the algorithm has concluded that the new direction is along the negative $x$-direction. If we do not begin at $(s',t')$, the program will just back-track along already found points of $S_k$.

It should be noted that while the algorithm has been presented in the two-dimensional case, it is easily generalized to apply to equations on $\R^{n}$. For instance, if we were to consider the case $f(x,y,z) = 0$, then a simple modification of Definition 1 including half-spheres, instead of half-balls, would allow for this algorithm to work with space curves.

\section{Verification}
\label{sec:astroid}
\begin{figure}
\begin{center}
\includegraphics[height=7cm] {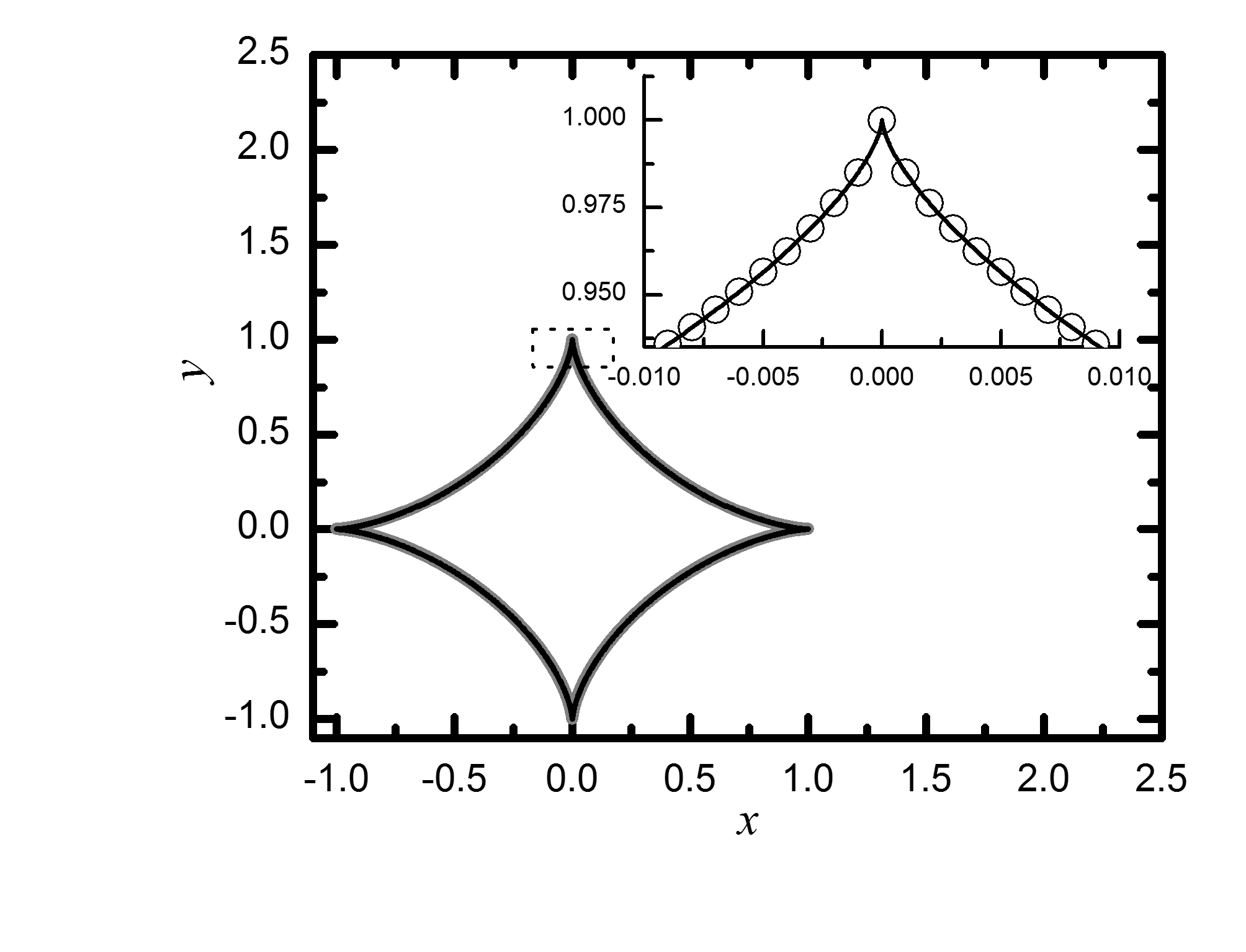}
\end{center}
\caption{Astroid plot generated using the turning point algorithm. Inset: parametrized curve (solid), approximated curve (circles)}
\label{fig:astroid}
\end{figure}

To demonstrate the effectiveness of the turning point algorithm, we chose to implement it on the astroid, a well known implicit equation given by $f(x,y) = x^{2/3}+y^{2/3}-1=0$. The resulting curve can be found in Figure \ref{fig:astroid}.

The accuracy of the curve generated with our algorithm was measured in percent error \eqref{PE}, and depended on several different parameters: $r$, the radius of $\partial R_{+}$, $k$ how many ``steps'' back the reference point was from the turning point, and $n$ the grid size of $\partial R_{+}$. Each parameter was varied until notable values were found. 

\begin{equation}\label{PE}
\PE(x_j,y_j) = \frac{\sqrt{x_j^2+y_j^2}-\sqrt{x^2+y^2}}{\sqrt{x^2+y^2}}\times 100\%
\end{equation}

It was found that a radius $r \in [0.0001\delta,100\delta]$ (where $\delta$ is the magnitude of the step-size $\Delta_x$ or $\Delta_y$), $1 \leq k \leq 10$ and $4 \leq n \leq 10$ yielded a maximum percent error $\ll 0.1\%$. However, for $n<4$ the algorithm failed to correctly navigate the turning point, and the plot ceased to continue.

\section{Example: Lubrication Model}
\label{sec:ode}

Consider the problem of determining the thickness of a thin film of lubricant inside a rotating cylinder. Using the lubrication approximation to the Navier-Stokes equations developed in \cite{moffatt, johnson}, one may describe the steady-state thickness, $h(\theta)$, by 
\begin{equation}
\label{eq:ode}
\frac{\epsilon}{3}(h'+h''') - \frac{1}{3} \cos(\theta) =  \frac{Q}{h^3} - \frac{1}{h^2}
\end{equation}
where $Q$ is the non-dimensional flux and $\epsilon$ is dependent on the surface tension and rate of rotation.
Coupling \eqref{eq:ode} with 
\begin{equation}
M = \int_{0}^{2\pi} h(\theta) \, d\theta
\label{eq:mass}
\end{equation}
it is not unreasonable to ask how the curve of points satisfying \eqref{eq:ode} and \eqref{eq:mass} behaves for a particular $\epsilon$. Because this curve may be characterized by a solution set to some equation $f_{\epsilon}(Q,M) = 0$, we will let $f_{\epsilon}$ denote the curve.

Since our intention is to demonstrate the robustness of the algorithm, and not the challenge behind solving $f_{\epsilon}(Q,M) = 0$ when either $Q$ or $M$ is fixed, we refer the curious reader to the works of Ashmore \emph{et{~}al.} \cite{ashmore} and Benilov \emph{et{~}al.} \cite{benilov1} for detailed exposition of that task. 

Briefly, we discretized $h$ on a uniform grid, and use Fourier spectral differentiation matrices which account for the $2\pi$-periodicity of $h$. A MATLAB program was written to accept very small values of $\epsilon$ and perform an iterated Newton-Raphson procedure to generate $f_\epsilon$ given a particular parameter and direction. Once a turning point was identified, the algorithm was performed.

\begin{figure}[ht!]
\centering
\includegraphics[height=7cm] {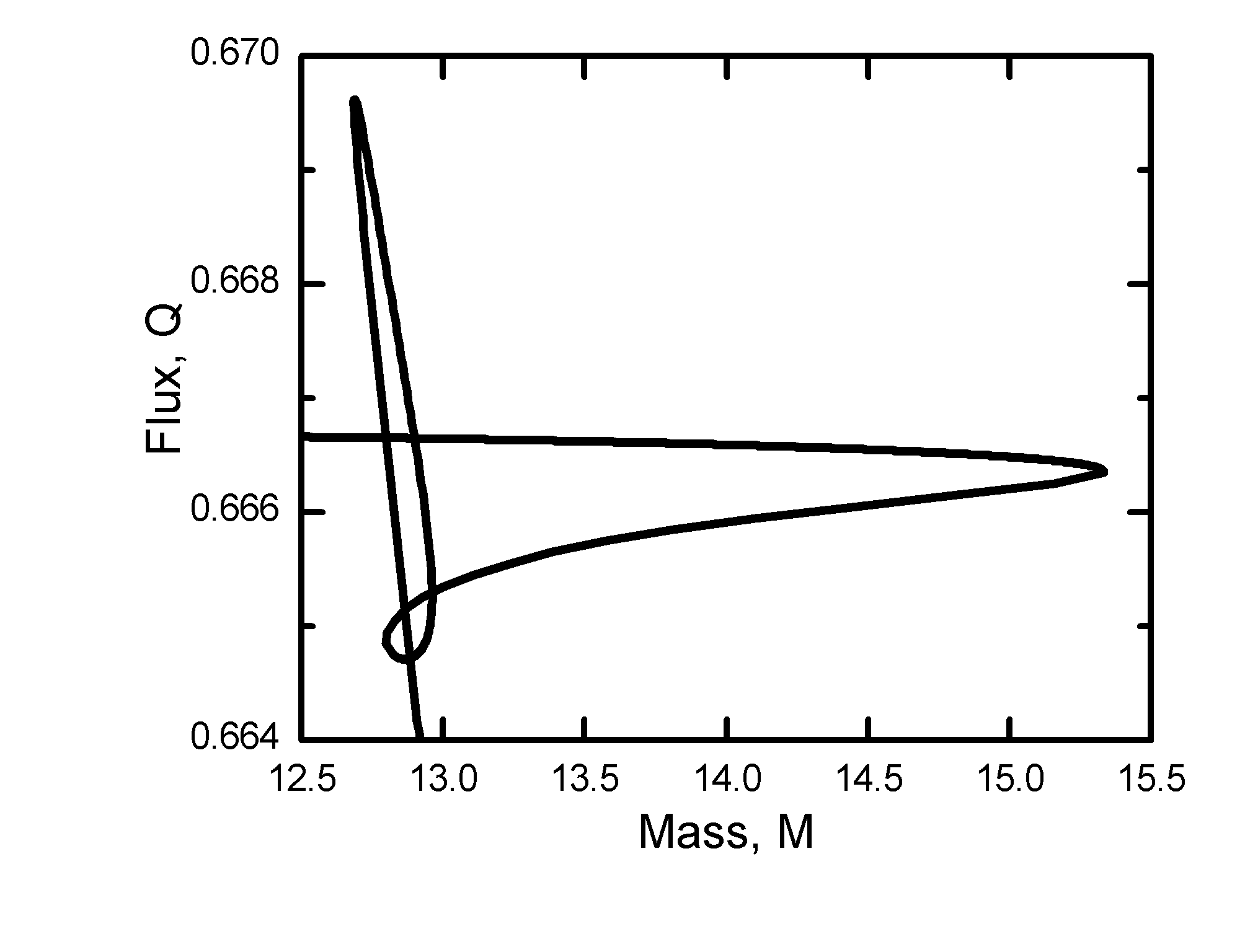}
\caption{A portion of the $(Q,M)$-bifurcation diagram generated for $\epsilon = 10^{-4}$. Modified from \cite{badali}.}
\label{fig:fluid}
\end{figure}

It was found that for $\epsilon\leq 10^{-3}$ complicated loop behavior similar to that exhibited in Figure \ref{fig:fluid} appeared. With the algorithm deployed to automate plotting, each successive turning point was over come with no difficulty. Interested readers may also refer to \cite{benilov1,benilov2,ashmore} for recent literature on this topic, as well as \cite{experiment1,experiment2,experiment3}. 

\section{Run-Time Analysis}
Due to the various settings in which the algorithm may be applied, a general run-time analysis is impossible. Instead, we will analyze the case when $f$ is not explicitly given, but inferred from differential or integral equations. In particular, we will perform a run-time analysis on the Lubrication Model example. As demonstrated in Section 5, it is clear that these settings are very useful in the study of certain ordinary and partial differential equations.

In this particular instance, we had two particular equations \eqref{eq:ode} and \eqref{eq:mass} which allowed us to solve for $Q$ and $M$ given a function $h : [0,2\pi] \to \R^{+}$. By approximating $h$ with $h_m$ (an $m\times 1$ column vector), and taking the appropriate Taylor expansion of elements of \eqref{eq:ode}, we were able to perform parameter continuation along the $M$-axis with a $m\times m$ matrix approximation of the operator $f$. 

With similar linearization techniques we were able to perform parameter continuation along the $Q$-axis with a $(m+1) \times (m+1)$ matrix approximation of $f$. (For more details about these techniques see Benilov \cite{benilov1}.)

In this setting, performing the Newton-Raphson method prescribed by Ashmore \cite{ashmore} and/or Benilov \cite{benilov1}, and thus completing Step 1 in the algorithm, is bottle-necked by the speed of matrix multiplication. 

Now, an algorithm for matrix multiplication of $m\times m$ matrices has been found with run-time $\mathcal{O}(m^{2.376})$ (i.e. the Coppersmith-Winograd algorithm\cite{cwalgo}). However, this algorithm is particularly hard to implement, and is usually used for theoretic bounds above all else. Thus, in the spirit of creating a practical run-time argument, we assume that Strassen's algorithm is used, and therefore the calculations in Step 1 can be carried out in $\mathcal{O}(m^{2.807})$ \cite{cormen}. Also, because Newton's method is usually performed in a for-loop for a fixed number of iterations, $\ell$ say, Step 1 may be preformed in $\mathcal{O}( n \,\ell \, m^{2.807}) = \mathcal{O}(n \,m^{2.807})$, where $n$ is the number of sample points on $\partial R_+$.  

Now, because calculating $\phi$ is reasonably assumed to take constant time, it is clear that the most cumbersome part of Step 2 is determining the point which minimizes $\phi$ (which is a sorting problem). Thus, if we assume $\r$ has $n$ points in it (the maximal amount), then Step 2 may be performed in $\Theta(n \log n)$. This follows from the fact that \cite{cormen} shows that that any comparison sort algorithm requires $\Omega(n \log n)$ in the worst case. Hence, taking any version of heap or merge-sort to perform Step 2 will require $\Theta(n \log n)$.

While there are a myriad of ways one may perform Step 3, any competent implementation should run in constant time since this step has only manipulates two points in $\R^2$. 

Hence, the algorithm can be found to run in $\mathcal{O}(n m^{2.807} + n \log n)$. However, as we saw in the astroid example, $n$ need not exceed 10, or so, points. Thus, if we approximate $h$ with any reasonable amount of resolution (i.e. $m = 2^7$ or $m = 2^8$), then $n \ll m$, and we may consider it something of a constant. Hence, the algorithm can be expected to run in  $\mathcal{O}(m^{2.807})$. 

\section{Conclusion}

It is easy to see that turning points make automated plotting a challenge, if not impossible. To overcome this issue we have outlined an algorithm that determines the direction iteration should continue in after reaching a turning point. As demonstrated with the astroid and thin-film differential equation examples, the algorithm overcomes turning points accurately, and effectively, while asking little in terms of computational power. It should be said that implicit functions and dynamical systems far from exhaust the areas where this algorithm may be helpful. With the amount of research being done today throughout the various sciences (pure and applied alike), it is clear that this algorithm can be applied anywhere data visualization is possible.

{\bf Acknowledgments}: This work was partially carried out during the Fields-MITACS Undergraduate Summer Research Program in 2010. The authors thank M. Chugunova for helpful discussions during this program.


\begin{thebibliography}{20}

\bibitem{moore1}
G. Moore and A. Spence,
\newblock {The calculation of turning points of nonlinear equations},
\newblock \emph{SIAM J. Numer. Anal.}
\newblock 17 (1980), pp. 567-576

\bibitem{moore2}
G. Moore and A. Spence,
\newblock {The convergence of operator approximations at turning points},
\newblock \emph{IMA J. Numer. Anal.}
\newblock 1 (1981), pp. 23-38

\bibitem{spence}
A. Spence and B. Werner,
\newblock {Non-simple turning points and cusps},
\newblock \emph{IMA J. Numer. Anal.}
\newblock 2 (1982), pp. 413-427

\bibitem{text}
R. Seydel,
\newblock \emph{Practical Bifurcation and Stability Analysis: From Equilibrium to Chaos},
\newblock Springer-Verlag, New York,
\newblock 1994

\bibitem{newton1}
E.J. Doedel, H.B. Keller, and J.P. Kernevez,
\newblock {Numerical analysis and control of bifurcation problems (I): Bifurcation in finite dimensions},
\newblock \emph{Internat. J. Bifur. Chaos Appl. Sci. Engrg.}
\newblock 1 (1991), pp. 493-520

\bibitem{newton2}
D.F. Davidenko,
\newblock {On a new method of numerical solution of systems of nonlinear equations},
\newblock \emph{MR}
\newblock 14 (1953), pp. 906

\bibitem{keller1}
H.B. Keller,
\newblock {Numerical solution of bifurcation and nonlinear eigenvalue problems},
\newblock in \emph{Applications of Bifurcation Theory},
\newblock Academic Press, New York,
\newblock 1977

\bibitem{keller2}
H.B. Keller,
\newblock \emph{Lectures on Numerical Methods in Bifurcation Problems},
\newblock Springer-Verlag, New York,
\newblock 1987

\bibitem{cliffe}
K.A. Cliffe, A. Spence, and S.J. Tavener,
\newblock {The numerical analysis of bifurcation problems with application to fluid mechanics},
\newblock \emph{Acta Numer.}
\newblock (2008), pp. 39-131

\bibitem{moffatt}
H.K. Moffatt,
\newblock {Behavior of a viscous film on the outer surface of a rotating cylinder},
\newblock \emph{J. de M\'{e}c.},
\newblock 16 (1977), pp. 651-673

\bibitem{johnson}
R.E. Johnson,
\newblock {Coating flow stability in rotating molding},
\newblock \emph{Engineering Science, Fluid Dynamics: A Symposium to Honor T.Y. Wu (ed. G. T. Yates)},
\newblock World Scientific,
\newblock 1990

\bibitem{benilov1}
E.S. Benilov, M.S. Benilov, and N. Kopteva,
\newblock {Steady rimming flows with surface tension},
\newblock \emph{J. Fluid Mech.}
\newblock 597 (2008), pp. 91-118

\bibitem{benilov2}
E.S. Benilov, M.S. Benilov, and S.B.G. O'Brian,
\newblock {Existence and stability of regularized shock solutions, with applications to rimming flows},
\newblock \emph{J. Eng. Math.}
\newblock 63 (2009), pp. 197-213

\bibitem{ashmore}
J. Ashmore, A.E. Hosoi, and H.A. Stone,
\newblock {The effect of surface tension on rimming flows in a partially filled rotating cylinder},
\newblock \emph{J. Fluid Mech.}
\newblock 479 (2003), pp. 65-98

\bibitem{experiment1}
R.E. Johnson,
\newblock {Steady-state coating flows inside a rotating horizontal cylinder},
\newblock \emph{J. Fluid Mech.}
\newblock 190 (1988), pp. 321-342

\bibitem{experiment2}
S.T. Thoroddsen and L. Mahadevan,
\newblock {Experimental study of coating flows in a partially-filled horizontally rotating cylinder},
\newblock \emph{Exp. Fluids}
\newblock 23 (1997), pp. 1-13

\bibitem{experiment3}
P.L. Evans, L.W. Schwartz, and R.V. Roy,
\newblock {Three-dimensional solutions for coating flow on a rotating horizontal cylinder: Theory and experiment},
\newblock \emph{Phys. Fluids}
\newblock 17 (2005), pp. 172102:1-20

\bibitem{badali}
D. Badali, M. Chugunova, D.E. Pelinovsky, and S. Pollack,
\newblock {Asymptotic behavior of regularized shock solutions in coating flows},
\newblock submitted to \emph{Phys. Fluids}
\newblock (January 2011)

\bibitem{cwalgo}
D. Coppersmith and S. Winograd,
\newblock {Matrix multiplication via arithmetic progressions},
\newblock \emph{J. Symbolic Comput.}
\newblock 9 (1990), pp. 251-280

\bibitem{cormen}
T.H. Cormen, C.E. Leiserson, R.L. Rivest, and C. Stein,
\newblock \emph{Introduction to Algorithms, Second Edition},
\newblock MIT Press and McGraw-Hill, Cambridge,
\newblock 2001

\end{thebibliography}
\end{document}